\documentclass[12pt]{elsart}
\usepackage{graphics}
\usepackage{graphicx}
\usepackage[dvips]{epsfig}
\usepackage{subfigure}
\usepackage{color}
\usepackage{amssymb,amsmath}

\newcommand{\R}{\mathbb{R}}

\newcommand{\RR}{\mathbb{R}}

\newenvironment{proof}{\noindent{\it Proof}\rm.}{\hfill $\Box$}
\newenvironment{proofof}[1]{\bigskip\noindent{\it Proof of~#1.}\rm}{\hfill $\Box$}

\theoremstyle{plain}
\newtheorem{theorem}{Theorem}[section]
\newtheorem{lemma}{Lemma} [section]
\newtheorem{corollary}[theorem]{Corollary}

\theoremstyle{definition}
\newtheorem{remarks}[theorem]{Remarks}
\newtheorem{remark}[theorem]{Remark}
\numberwithin{equation}{section}
\begin{document}
\begin{frontmatter}

\title{Differential inequalities for Riesz means and Weyl-type bounds for eigenvalues\thanksref{gracias}}
\author{Evans M. Harrell II}
\ead{harrell@math.gatech.edu}
\address{School of Mathematics, Georgia Institute of
Technology, Atlanta, GA 30332-0160 U.S.A.}
\author{Lotfi Hermi} \ead{hermi@math.arizona.edu}
\address{Department of Mathematics,
University of Arizona, 617 Santa Rita, Tucson, AZ 85721 U.S.A.}
\thanks[gracias]{The second author would like to thank the
Georgia Tech School of Mathematics for hospitality and support
while doing this work.}

\begin{abstract}
We derive differential inequalities and difference inequalities
for Riesz means of eigenvalues of the Dirichlet Laplacian,
$$R_{\sigma}(z) := \sum_k{(z - \lambda_k)_+^{\sigma}}.$$
Here $\left\{ \lambda_k\right\}_{k=1}^{\infty}$
are the ordered eigenvalues of the Laplacian on a bounded domain
$\Omega \subset \R^d$,  and $x_+ := \max(0, x)$
denotes the positive part of the quantity $x$.  As corollaries of
these inequalities, we derive Weyl-type bounds on $\lambda_k$, on
averages such as $\overline{\lambda_k} := {\frac 1 k }\sum_{\ell
\le k}\lambda_\ell$, and on the eigenvalue counting function.  For
example, we prove that for all domains and all $k \ge j
\frac{1+\frac d 2}{1+\frac d 4}$,
$$
{\overline{\lambda_{k}}}/{\overline{\lambda_{j}}} \le
2 \left(\frac{1+\frac d 4}{1+\frac d 2}\right)^{1+\frac 2 d}\left( {\frac k j}\right)^{\frac 2 d}.
$$
\end{abstract}

\end{frontmatter}

\section{Introduction} \label{intro}

This article is concerned with the spectrum of the Laplace
operator $-\Delta$ on a bounded Euclidean domain
$\Omega \subset \R^d$, $d\ge 1$, with vanishing
Dirichlet boundary conditions on $\partial \Omega$, {\it i.e.},
the classic problem of the vibratory modes of a finite elastic
body with fixed boundary. With straightforward changes it is also
possible to treat Schr\"odinger operators $- \Delta + V(x)$ with
discrete spectra.  Indeed,
virtually no changes at all are necessary when $V(x) \ge 0$,
other than possibly the upper bound in
\eqref{RUprLower}.
For simplicity of exposition, however,
we refer below only to the Laplacian.

The spectrum $\left\{
\lambda_k\right\}_{k=1}^{\infty}$, of $-\Delta$ is subject to the
Weyl law, \cite{Weyl,CH,BH,lieb1}
\begin{equation} \label{Weyl}
\lambda_k \sim 4 \pi \frac{\left(\Gamma(1 +
\frac{d}{2})\right)^{\frac 2 d} }{|\Omega|^{\frac{2}{d}}} \,
k^{\frac{2}{d}}
 \end{equation}
as $k \to \infty$; to the Berezin-Li-Yau inequality,
\cite{berezin,LY,Lap1,LapWeidl1,LapWeidl2}
\begin{equation} \label{LiYau}
{\overline{\lambda_{k}}}  := {\frac 1 k }\sum_{\ell \le k}\lambda_\ell
\ge 4 \pi
\frac{\left(\Gamma(1 + \frac{d}{2})\right)^{\frac 2 d}
}{|\Omega|^{\frac{2}{d}}}  \, {\frac{k^{\frac 2 d}}{1+{\frac
2 d}}};
\end{equation}
and to other familiar constraints relating the spectrum to
the geometry of $\Omega$. The Berezin-Li-Yau inequality is a
notable example of an inequality for which both sides are of the
same order in $k$, when the asymptotic expression in \eqref{Weyl}
is substituted for $\lambda_k$.  A formula or expression where
this is the case is said to be of {\it Weyl-type}.

At the same time, the spectrum is subject to ``universal bounds''
by which certain expressions involving eigenvalues dominate others
with no reference to the geometry of $\Omega$.   The expressions
occurring in known universal bounds include moments of eigenvalues, in
particular sums $\sum_{j=1}^k{\lambda_j}$ as in Berezin-Li-Yau.
(See \cite{A1,Her} for a review of universal
spectral bounds.) Although early universal bounds like that of
Payne, P\'olya, and Weinberger \cite{PPW2} were not of Weyl-type,
universal bounds of Weyl-type have been known since H.~C. Yang's
unpublished 1991 article
\cite{Yang,HS,LevPar,A1,AH2,ChengYang1,ChengYang2,AH3,AH4,Her}.

One of the goals here is to show that trace identities of the type
introduced in \cite{HS} imply tight, Weyl-type bounds on ratios of
eigenvalues belonging to Laplace spectra. To our knowledge this
has not been much explored except by Hermi \cite{Her} and by Cheng
and Yang in \cite{ChengYang2}, with some remarks in
\cite{ChengYang1}.  Analogous bounds for Schr\"odinger operators
in one dimension can be found in the work of Ashbaugh and Benguria
\cite{AB89}.  Ratios involving averages ${\overline{\lambda_{k}}}$
will arise in a natural way through this analysis.

The earliest bound for the ratio of a large eigenvalue to the
fundamental eigenvalue is due to Ashbaugh and Benguria \cite{AB4},
who proved
\begin{equation}
\dfrac{\lambda_{2^m}}{\lambda_1} \le
\left(\dfrac{j_{d/2,1}^2}{j_{d/2-1,1}^2}\right)^m. \label{ab94}
\end{equation}
While optimal for low-lying eigenvalues, \eqref{ab94} is not of
Weyl-type since the right side of the inequality
behaves like $k^{5.77078/d}$ as $k\to \infty$ (see the details in
\cite{Her}). Weyl-type bounds were proved in \cite{Her},  of the
form
\begin{equation}
\frac{\lambda_{k+1}}{\lambda_1}\le 1 +
\left(1+\dfrac{d}{2}\right)^{2/d} H_d^{2/d} \, k^{2/d},
\label{her1}
\end{equation}
and
\begin{equation}
\frac{\overline{\lambda}_k}{\lambda_1} \le 1 + \frac{H_d^{2/d}}
{1+\frac{2}{d}} \, k^{2/d}, \label{her2}
\end{equation}
where
\begin{equation} \label{chiti-constant}
H_d=\dfrac{2 \ d}{j_{d/2-1,1}^2
J_{d/2}^2(j_{d/2-1,1})}.\end{equation}
As usual, $j_{\alpha,p}$ denotes
the $p$-th positive zero of the Bessel function
$J_{\alpha}(x)$ \cite{AS}.   In \cite{ChengYang2} Cheng and Yang
prove that
 \begin{equation} \label{CYbound}
\frac{\lambda_{k+1}}{\lambda_1} \le \left(1 + \frac{4}{d}\right) k^{\frac 2 d},
\end{equation}
as well as some incremental improvements for large values of $k,d$ of the form
 \begin{equation}\label{CYcomplicated}
\frac{\lambda_{k+1}}{\lambda_1} \le C_0(d,k) k^{\frac 2 d}.
\end{equation}
For instance, for $k\ge d+1$, Cheng and Yang
\cite{ChengYang2} improve \eqref{CYbound} to
\begin{equation}\label{CYbound2}
\frac{\lambda_{k+1}}{\lambda_1} \le \left(1 + \frac{4}{d}\right)
\, \left(1+ \frac{8}{d+1}+ \frac{8}{(d+1)^2}\right)^{\frac 1 2} \,
\left(d+1\right)^{-\frac 2 d} \, k^{\frac 2 d}.
\end{equation}

For comparison, by combining the Weyl law with the
Rayleigh-Faber-Krahn inequality \cite{Faber,Krahn1,Krahn2}, there
immediately results an asymptotic upper bound on $\lambda_k /
\lambda_1$, of the form
\begin{equation} \label{FKbound} \frac{4 \pi}{\lambda_1^*} \,
\frac{\left(\Gamma(1 + \frac{d}{2})\right)^{\frac 2 d}
}{|\Omega|^{\frac{2}{d}}} \, k^{\frac{2}{d}} = \frac{4
\left(\Gamma(1 + \frac{d}{2})\right)^{\frac 4
d}}{j_{\frac{d}{2}-1,1}^2} k^{\frac{2}{d}},
\end{equation}
where
\[\lambda_1^* =
\frac{\pi \, j_{\frac{d}{2}-1,1}^2}{\left(\Gamma(1 + \frac{d}{2})
\, |\Omega|\right)^{\frac{2}{d}}} \] is the explicit value of the
fundamental eigenvalue when $\Omega$ is a ball.

Ideally, a Weyl-type bound would contain a constant commensurate
with that on the right side of \eqref{FKbound}. The bound of
\eqref{CYbound} is numerically nearly 4.34 times as large as the
ideal when $d=2$. It will be shown below that the constants in
\eqref{CYbound} and \eqref{CYbound2} can be reduced. (See Table
\ref{table2} below.)

Our technique will make use of differential inequalities and
difference inequalities for Riesz means of eigenvalues. Safarov,
Laptev, and Weidl have long advocated Riesz means as a tool for
understanding inequalities like those of Lieb-Thirring and
Berezin-Li-Yau, and we draw some of our inspiration from
\cite{Saf1,Saf2,Lap1,Lap2,LapWeidl2,LapWeidl1}.

\medskip

\section{Differential and difference inequalities for Riesz means of eigenvalues} \label{body}

For background information we refer to the monograph of
Chandrasekharan and Minakshisundaram \cite{ChandraMinak}, where
Riesz means are referred to as {\it typical means}. Recall that if
$\{\lambda_k\}_{k=1}^{\infty}$ is an increasing sequence of real
numbers, then for any real $z$, the {\it Riesz
mean of order $\sigma$} of $\{\lambda_k\}$ can be defined as
 \begin{equation} \label{Riesz}
R_{\sigma}(z) := \sum_k{(z - \lambda_k)_+^{\sigma}}
 \end{equation}
for $\sigma > 0$.
When $\sigma = 0$, we interpret $R_0(z)$ as the eigenvalue-counting function
$R_0(z) := {\mathcal N}(z) := \lim_{\sigma \downarrow 0}{R_{\sigma}(z)}$.

One of the main results consists of differential inequalities for
$R_{\sigma}(z)$ with respect to $z$ and difference inequalities
for $R_{\sigma}(z)$ with respect to $\sigma$.

\begin {theorem} \label{BasicIneq}
For $0 < \sigma \le 2$ and $z \ge \lambda_1$,
 \begin{equation} \label{Difference1}
 R_{\sigma-1}(z) \ge \left({1 + {\frac d 4}}\right) {\frac 1 z} R_{\sigma}(z);
 \end{equation}
  \begin{equation} \label{Differential1}
 R_{\sigma}^{\prime}(z) \ge \left({1 + {\frac d 4}}\right) {\frac \sigma z} R_{\sigma(z)};
 \end{equation}
and consequently \[\dfrac{R_{\sigma}(z)}{
z^{\sigma+{\frac {d \sigma} 4}}}
\]
is a nondecreasing function of $z$.
\medskip
\par\noindent
For $2 \le \sigma < \infty$ and $z \ge \lambda_1$,
 \begin{equation} \label{Difference2}
 R_{\sigma-1}(z) \ge \left(1 + {\frac d {2 \sigma}}\right) {\frac 1 z} R_{\sigma}(z);
 \end{equation}
  \begin{equation} \label{Differential2}
 R_{\sigma}^{\prime}(z) \ge  \left(\sigma + {\frac d {2}}\right) {\frac 1 z} R_{\sigma}(z);
 \end{equation}
and consequently \[\dfrac{R_{\sigma}(z)}
{z^{\sigma+{\frac d 2}}} \] is a nondecreasing function of $z$.
\end{theorem}

\begin{remarks}\label{rk1}
\rm\item{1.}
The values of $z$ are restricted for the simple reason that
$R_{\sigma}(z) = 0$ for $\sigma > 0, z \le \lambda_1$, and is undefined for
$\sigma \le 0, z \le \lambda_1$.
\item{2.}
It is only necessary to prove \eqref{Difference1} and \eqref{Difference2}, since
from \eqref{Riesz},
$$
 R_{\sigma}^{\prime}(z) = \sigma R_{\sigma-1}(z),
$$
and since solving the elementary differential inequalities
\eqref{Differential1} and \eqref{Differential2} easily yields the
remaining statements about functions of the form
$R_{\sigma}(z)/z^p$. \item{3.}  In the case $\sigma=2$,
\eqref{Difference2} is related to the inequality of H.~C. Yang
\cite{Yang,HS,LevPar,AH3,HaHer1}. \item{4.} Differential
inequalities were considered in \cite{HS} for the partition
function (or trace of the heat kernel) $Z(t) := {\tt tr}\left(e^{t
\Delta}\right)$, where it was shown that $Z(t) t^{d/2}$ is a
nonincreasing function. In a future article \cite{HaHer1}, it will
be shown that the differential inequalities of Theorem
\ref{BasicIneq} imply the equivalence of the Berezin-Li-Yau
inequality and a classical inequality due to Kac \cite{Kac} about
the partition function.

\end{remarks}
To prepare the proof we state two lemmas.
\begin{lemma}\label{HS}
Denoting the $L^2$-normalized eigenfunctions of the Laplace operator
$\left\{u_j\right\}$, let
$$T_{\alpha jm} := \left|{\left({\frac{\partial u_j}{\partial x_{\alpha}}, u_m}\right)}\right|^2$$
for $j, m = 1, \dots$ and $\alpha = 1, \dots, d$.  Then for each fixed $\alpha$,
\begin{equation}\label{HSIneq}
R_\sigma(z)
= 2\sum\limits_{j,m: \lambda_j\ne \lambda_m}
{{\frac{\left(z - \lambda_j\right)_+^\sigma - \left(z - \lambda_m\right)_+^\sigma}
{\lambda_m-\lambda _j}}\,
{T}_{\alpha jm}}+4\sum\limits_{j,q: \lambda _j \le z <
\lambda_q}
{{\frac{\left(z -\lambda _j\right)^\sigma}
{{\lambda }_{q}-{\lambda }_{j}}}\,{T}_{\alpha jq}}.
\end{equation}
\end{lemma}

This lemma is the trace identity of Harrell and Stubbe (\cite{HS}
Theorem 1, (4)), specialized to $f(\lambda) = (z -
\lambda)_+^\sigma$.   Versions of the trace identity of \cite{HS}
also appear in some later articles, e.g.,
\cite{LevPar,AH2,AH3,EHI}.

\begin{lemma}\label{Elem}
Let $0 < x < y$ and $\sigma \ge 0$.  Then
\begin{equation}
\frac{y^\sigma - x^\sigma}{y-x} \le C_{\sigma}\left({y^{\sigma-1} + x^{\sigma-1}}\right),
\end{equation}
where
\begin{equation*}
C_{\sigma} :=
   \begin{cases}
   \frac{\sigma}{2},      &\text{if $0 \le \sigma < 1$}\\
   1,     &\text{if $1 \le \sigma \le 2$}\\
   \frac{\sigma}{2},      &\text{if $2 \le \sigma < \infty$.}
   \end{cases}
\end{equation*}
\end{lemma}

\noindent
\begin{proof}
By a scaling, it suffices to assume $x=1$.  We then seek the supremum of
$$\frac{y^{\sigma} - 1}{(y-1)(y^{\sigma-1} + 1)}$$
for $1 < y < \infty$.  A calculus exercise shows that the supremum is approached as
$y \downarrow 1$ when $0 < \sigma < 1$ or $2 \le \sigma < \infty$, whereas
it is approached as $y \to \infty$ when $1 \le \sigma \le 2$.  (It's constant
when $\sigma=2$.)  The stated values are
obtained with l'H\^opital's rule in the cases $0 < \sigma < 1$ and $2 \le \sigma < \infty$.
\end{proof}
\medskip

\begin{proofof}{Theorem~\ref{BasicIneq}}
Let the first term on the right of \eqref{HSIneq}
be
\[
G(\sigma,z, \alpha):=2\sum\limits_{j,m: \lambda_j\ne \lambda_m}
{{\frac{\left(z - \lambda_j\right)_+^\sigma - \left(z -
\lambda_m\right)_+^\sigma} {\lambda_m-\lambda _j}}\, {T}_{\alpha
jm}}.
\]
By Lemma~\ref{Elem}, this expression simplifies to
\begin{eqnarray*}
G(\sigma,z, \alpha)&&\quad = 2 \sum\limits_{j,m: \lambda_{j,m} \le
z, \lambda_j\ne \lambda_m} {{\frac{\left(z -
\lambda_j\right)^\sigma - \left(z - \lambda_m\right)^\sigma} {(z -
\lambda_j ) - (z - \lambda_m)}}\, {T}_{\alpha jm}} \\
&&\quad\le 2 \, C_\sigma \sum\limits_{j,m: \lambda_{j,m} \le z}
{{\left((z - \lambda_j)^{\sigma-1} + (z -
\lambda_m)^{\sigma-1}\right)} \,
{T}_{\alpha jm}} \\
&&\quad = 4 \, C_\sigma \sum\limits_{j,m: \lambda_{j,m} \le z}
{{{\left(z - \lambda_j\right)_+^{\sigma-1}} }\, {T}_{\alpha jm}}
\end{eqnarray*}
by symmetry in $j \leftrightarrow m$.  Extending the sum to all
$m$ and subtracting the same quantity from the final term in
\eqref{HSIneq}, we find
\begin{equation}\label{P1}
R_\sigma(z) \le 4 C_\sigma \sum\limits_{j: \lambda_{j} \le z,\; {\tt
all }\; m} {{{\left(z - \lambda_j\right)_+^{\sigma-1}} }\,
{T}_{\alpha jm}} + 4 H(\sigma,z, \alpha),
\end{equation}
where
\begin{equation}\label{P2}
H(\sigma,z, \alpha) :=
\sum\limits_{j,q: \lambda _j \le z <\lambda_q}{{T}_{\alpha jq} \left(z -\lambda _j\right)^{\sigma-1} \left({\frac{\left(z -\lambda _j\right)-C_\sigma \left(\lambda_q -\lambda _j\right)}
{{\lambda }_{q}-{\lambda }_{j}}}\right)}.
\end{equation}

Next observe that because $\left\{u_m\right\}$ is a complete orthonormal set,
$$\sum_m{T_{\alpha jm}} = \left\|{\frac{\partial u_j}{\partial x_{\alpha}}}\right\|_2^2$$
and thus
$$\sum_{m,\alpha} {T_{\alpha jm}} = \|\nabla u_j\|_2^2 = \lambda_j.$$
Therefore we may average over $\alpha = 1, \dots, d$ in \eqref{P1} to obtain
$$
R_\sigma(z) \le \frac{4 C_\sigma}{d} \sum\limits_{j}
{{{\left(z - \lambda_j\right)_+^{\sigma-1}}
}\,
\lambda_j} +  \frac{4}{d} \sum_{\alpha=1}^d{H(\sigma,z, \alpha)},
$$
or, since
$$
 \sum\limits_{j}
{{{\left(z - \lambda_j\right)_+^{\sigma-1}}
}\,\lambda_j} =
z R_{\sigma-1}(z) - R_{\sigma}(z),
$$
\begin{equation}\label{P3}
\left(1 + \frac{4 C_\sigma}{d}\right) R_\sigma -
\frac{4 z C_\sigma}{d} R_{\sigma-1}(z)
\le
\frac{4}{d} \sum_{\alpha=1}^dH(\sigma,z, \alpha).
\end{equation}

\par\noindent
We consider three cases:

\par \noindent
Case I.  $1 \le \sigma \le 2$.  In this case $C_\sigma = 1$ and $H(\sigma,z, \alpha) \le 0$,
establishing that
$$
\left(1 + \frac{4}{d}\right) R_\sigma - \frac{4 z }{d} R_{\sigma-1}(z) \le 0,
$$
which is equivalent to \eqref{Difference1} for this range of $\sigma$.

\par \noindent
Case II.  $0 < \sigma < 1$.

Since the sum defining $H$ runs over $\lambda_q > z$,
$$
\frac{\left(z -\lambda _j\right)-C_\sigma \left(\lambda_q -\lambda
_j\right)} {\lambda_q-\lambda_j} \le \frac{\left(\lambda_q
-\lambda _j\right)-C_\sigma \left(\lambda_q -\lambda _j\right)}
{\lambda _q-\lambda_j} = 1 - C_\sigma.
$$

Therefore
$$
H(\sigma,z, \alpha) \le (1 - C_\sigma)
\sum\limits_{j,q: \lambda _q \ge z}
{{T}_{\alpha jq} \left(z -\lambda _j\right)_+^{\sigma-1} },
$$
and since in this case $1 - C_\sigma = 1 - \sigma/2 > 0$,
we may extend the sum over all $q$,
obtaining
$$
\sum_{\alpha=1}^d{H(\sigma,z, \alpha) \le (1 - C_\sigma)
\sum\limits_{j}
{\lambda_j} \left(z -\lambda _j\right)_+^{\sigma-1} }.
$$
Substituting this into \eqref{P3}, there is a cancellation of
the $C_\sigma$'s, and we again obtain
$$
\left(1 + \frac{4}{d}\right) R_\sigma - \frac{4 z }{d} R_{\sigma-1}(z) \le 0,
$$
equivalent to \eqref{Difference1}.

\par \noindent
Case III.  $\sigma > 2$.
Arguing as in Case II, we come as far as
$$
H(\sigma,z, \alpha) \le (1 - C_\sigma)
\sum\limits_{j,q: \lambda _j \le z <
\lambda_q}
{{T}_{\alpha jq} \left(z -\lambda _j\right)^{\sigma-1} },
$$
but since now $(1 - C_\sigma) = 1 - \sigma/2 < 0$,
we cannot extend the sum in $q$ to simplify
${T}_{\alpha jq}$, and only conclude that
$$
\left(1 + \frac{4 C_\sigma}{d}\right) R_\sigma - \frac{4 z C_\sigma}{d} R_{\sigma-1}(z)
\le 0,
$$
which is \eqref{Difference2}.

\end{proofof}

Convexity provides insight into why the proof divides into three
ranges of $\sigma$. The Riesz mean $R_\sigma(z)$ loses convexity
in $z$ when $\sigma < 1$, and its derivative loses convexity
already when $\sigma < 2$.  Thus both $\sigma = 1$ and $2$ are
values at which necessary inequalities in the proof reverse.

Since by the theorem, $R_\sigma(z) z^{-p}$ is a nondecreasing
function for appropriate values of $p$, we obtain
lower bounds of the form $R_\sigma(z) \ge C z^p$ for all $z \ge
z_0$ as soon as $R_\sigma(z_0)$ is known, or estimated from below.
Upper bounds can be obtained from the limiting behavior of
$R_\sigma(z)$ as $z \to \infty$, as given by the Weyl law.

\begin{corollary}\label{Rbds}
For all $\sigma \ge 2$ and $z \ge \left(1 + {\frac{2 \sigma}{d}}\right)\lambda_1$,
\begin{equation}\label{RUprLower}
\left(\frac{2 \sigma}{d}\right)^\sigma
\lambda_1^{-\frac{d}{2}}\left(\frac{z}{1+\frac{2
\sigma}{d}}\right)^{\sigma + \frac{d}{2}} \le R_\sigma(z) \le
L_{\sigma, d}^{\it cl} |\Omega| z^{\sigma + \frac{d}{2}},
\end{equation}
where
 \begin{equation} \label{ClassCst}
L_{\sigma, d}^{\it cl} := \frac{\Gamma(\sigma+1)}{(4 \pi)^{\frac{d}{2}} \Gamma\left(\sigma+1+\frac{d}{2}\right)}.
\end{equation}
\end{corollary}

\begin{remark}\label{rk2} \rm
The upper bound on $R_\sigma(z)$ is not new, but is a result of
Laptev and Weidl \cite{Lap1,LapWeidl2}, for which we provide an
independent method of proof for $\sigma \ge 2$. They regard it as
a version of the Berezin-Li-Yau inequality, to which it is
directly related by a Legendre transform when $\sigma=1$. See also
\cite{HaHer1} for a discussion of the various connections.
\end{remark}

\begin{proof}
Since $R_{\sigma}(z) \ge (z_0 -
\lambda_1)_+^{\sigma}$, for any $z_0
> \lambda_1$, it follows from {Theorem~\ref{BasicIneq}} that for
all $z \ge z_0$,
$$
R_{\sigma}(z) \ge (z_0 - \lambda_1)_+^{\sigma}
\left(\frac{z}{z_0}\right)^{\sigma + \frac{d}{2}}.
$$
By a straightforward calculation, the coefficient of $z^{\sigma +
\frac{d}{2}}$ is maximized when $z_0 =  \left(1 + {\frac{2
\sigma}{d}}\right)\lambda_1$, and the lower bound
that results is the left side of \eqref{RUprLower}.

As for the other inequality, note, following \cite{Lap1,LapWeidl2}, that the Weyl law implies that
$$
\frac{R_{\sigma}(z)}{z^{\sigma + \frac{d}{2}}} \to L_{\sigma,
d}^{\it cl} \, |\Omega|
$$
as $z \to \infty$.
Since $\frac{R_{\sigma}(z)}{z^{\sigma + \frac{d}{2}}}$ is a
nondecreasing function, it is less than this limit for all finite
$z$.
\end{proof}

{\begin{remark}\label{rk-added} \rm It is also possible to prove
the upper bound in \eqref{RUprLower} by invoking the difference
inequality \eqref{Difference2} directly. Note that
\eqref{Difference2} can be rewritten in the form
\begin{equation}\label{Aizen-Lieb-Form}
\dfrac{R_{\sigma-1}(z)}{L_{\sigma-1,d}^{cl} \,
z^{\sigma-1+\frac{d}{2}}} \ge
\dfrac{R_{\sigma}(z)}{L_{\sigma,d}^{cl} \,
z^{\sigma+\frac{d}{2}}}, \end{equation}
by virtue of the fact that
\[L_{\sigma-1,d}^{cl}=\left(1+\frac{d}{2 \sigma}\right) \,
L_{\sigma,d}^{cl}.
\]
Ineq.~\eqref{Aizen-Lieb-Form} is a version of a monotonicity
principle of Aizenman and Lieb \cite{AS} at the level of the
ratios rather than their suprema. The proof is
then complete owing to the monotonicity principle for
$R_{\sigma}(z)/z^{\sigma+d/2}$ and the Weyl asymptotic law as
applied to $R_{\sigma-1}(z)/z^{\sigma-1+d/2}$, i.e., for $z\ge
z_0> 0$
\[\dfrac{R_{\sigma-1}(z)}{L_{\sigma-1,d}^{cl} \,
z^{\sigma-1+\frac{d}{2}}} \ge
\dfrac{R_{\sigma}(z)}{L_{\sigma,d}^{cl} \, z^{\sigma+\frac{d}{2}}}
\ge \dfrac{R_{\sigma}(z_0)}{L_{\sigma,d}^{cl} \,
z_0^{\sigma+\frac{d}{2}}}.
\]
Sending $z\to\infty$ results in the Berezin-Li-Yau inequality for
$\sigma\ge2$ once more.
\end{remark}

Actually, the monotonicity principle described in Theorem
\ref{BasicIneq} implies that for sufficiently large $z_0$, the
lower bound on $R_{\sigma}(z)/z^{\sigma+ \frac d 2}$ is
arbitrarily close to the upper bound.  That is, given any domain
$\Omega$ and any $\epsilon > 0$, there exists a large but finite
$z_{\epsilon}$ such that for $z \ge z_{\epsilon}$,
$$R_{\sigma}\ge \left(L_{\sigma, d}^{\it cl} - \epsilon\right) |\Omega| z^{\sigma + \frac{d}{2}}.$$
In this estimate, however, $z_{\epsilon}$ is not independent of
the shape of $\Omega$.

Arguments nearly identical to those adduced for
{Corollary~\ref{Rbds}} lead to lower bounds on $R_{\sigma}$ when
$\sigma < 2$, but they are not of Weyl-type due to the different
exponent of $z$ in {Theorem~\ref{BasicIneq}}. To remedy this
deficiency, we call on \eqref{Difference1} and  \eqref{Difference2}:

\begin{corollary}\label{Rlowerbd<2}
For $1 \le \sigma  < 2$ and $z \ge \left(1 + {\frac{2 \sigma + 2}{d}}\right)\lambda_1$,
\begin{equation}\label{P4}
R_\sigma(z) \ge
\frac{\left(2 \sigma + 2\right)^{\sigma} d^{\frac d 2}}{\left(d + 2 \sigma + 2\right)^{\sigma + \frac{d}{2}}}  \lambda_1^{-\frac{d}{2}} z^{\sigma+ \frac{d}{2}},
\end{equation}
and for $0 \le \sigma < 1$, $z\ge \left(1 + {\frac{2 \sigma + 4}{d}}\right)\lambda_1$,
\begin{equation}\label{P5}
R_\sigma(z) \ge
\frac{\left(1 + \frac{d}{4} \right)\left(2 \sigma + 4\right)^{\sigma+1} d^{\frac d 2}}{\left(d + 2 \sigma + 4\right)^{\sigma + 1+\frac{d}{2}}}  \lambda_1^{-\frac{d}{2}} z^{\sigma+ \frac{d}{2}}.
\end{equation}

\end{corollary}

For convenience we collect the most important cases $\sigma=0, 1$,
which, by Theorem \ref{BasicIneq}, simplify to:

\begin{equation}\label{P5.5}
R_1(z) \ge \left(1 + \frac{d}{4}\right) \frac{1}{z} R_2(z) \ge
\frac{4 \ d^{\frac d 2}}{\left(d + 4 \right)^{1 + \frac{d}{2}}}  \lambda_1^{-\frac{d}{2}} z^{1+ \frac{d}{2}},
\end{equation}
and,
\begin{equation}\label{P6}
{\mathcal N}(z) = R_0(z)  \ge \left(1 + \frac{d}{4}\right)^2 \frac{1}{z^2} R_2(z) \ge \left(\frac{z}{\left(1 + \frac{4}{d}\right) \lambda_1}\right)^{\frac d 2}.
\end{equation}

\begin{remark}\label{rk-added2} \rm
The Cheng-Yang bound \eqref{CYbound} is a simple corollary of
\eqref{P6}; for the proof let $z$ approach $\lambda_{k+1}$ and
collect terms.  These inequalities and the Cheng-Yang bound will
be improved and generalized below. Furthermore, the lower bounds
provided by \eqref{RUprLower} and \eqref{P4} compete with the
bound
\begin{equation}
R_{\sigma}(z) \ge H_d^{-1} \lambda_1^{-d/2} \
\dfrac{\Gamma(1+\sigma) \Gamma(1+d/2)}{\Gamma(1+ \sigma+ d/2)} \
\left(z-\lambda_1\right)_{+}^{\sigma+d/2}, \label{riesz-low}
\end{equation}
valid for $\sigma\ge 1$, where $H_d$ is defined by
\eqref{chiti-constant} above. The case $\sigma=1$ of
\eqref{riesz-low} can be found in \cite{Her} and is in fact hidden
in earlier work of Laptev \cite{Lap1}. It is the Legendre
transform of this inequality that provides the bound \eqref{her2}
which competes with Cheng and Yang's \eqref{CYbound}. Additional
discussion and comparisons appear in Section \ref{comparison}
(see also \cite{HaHer1}).
\end{remark}

One option for generalizing Inequalities \eqref{P6} and
\eqref{CYbound} is to compare with parts of the spectrum lying
above $\lambda_1$. To avoid complications, we shall focus on
inequalities involving $R_2(z)$. Just as $\overline{\lambda_j}$
denotes the mean of the eigenvalues $\lambda_\ell$ for $\ell \le
j$, define
$$
\overline{\lambda_j^2} := {\frac 1 j
}\sum_{\ell \le j}\lambda_{\ell}^2.
$$
Then setting $j = [z]$
(greatest integer $\le z$), we find that
\begin{equation}
R_2(z) = [z] \left(z^2 - 2 z {\overline{\lambda_{[z]}}} + \overline{\lambda_{[z]}^2}\right).
\end{equation}
For any integer $j$ and all $z \ge j$,
$$R_2(z)\ge Q(z,j) :=  j \left(z^2 - 2 z {\overline{\lambda_{j}}} + \overline{\lambda_{j}^2}\right).$$
Using {Theorem~\ref{BasicIneq}}, for $z \ge z_j \ge \lambda_j$,
\begin{equation}\label{P7}
R_2(z)\ge Q(z_j,j) \left(\frac{z}{z_j}\right)^{2+\frac{d}{2}}.
\end{equation}

A good -- simple but not optimized -- choice is $z_j = \left(1 +
{\frac 4 d}\right) {\overline{\lambda_{j}}}$. Observe that because
of the Cauchy-Schwarz inequality, ${\overline{\lambda_{j}}}^2 \le
{\overline{\lambda_{j}^2}}$, so
$$ Q(z,j) = j  \left( \left(z - {\overline{\lambda_{j}}} \right)^2 + \overline{\lambda_{j}^2} - \overline{\lambda_{j}}^2\right)
$$
\begin{equation}\label{P8}
\quad\quad \ge j \left(z - {\overline{\lambda_{j}}} \right)^2.
\end{equation}

This establishes the following.
\begin{corollary}\label{lambdabarbds}
Suppose that
$z \ge \left(1 + {\frac 4 d}\right) {\overline{\lambda_{j}}}$.  Then

\begin{equation}\label{Rbddbylambdabar}
R_2(z) \ge \frac{j z^{2+{\frac d 2}}}{\left(1+{\frac d 4}\right)^2 \left(\left(1+{\frac 4 d}\right) {\overline{\lambda_{j}}}\right) ^{\frac{d}{2}}},
\end{equation}
and therefore,
\begin{equation}\label{P9}
R_1(z) \ge \frac{j z^{1+{\frac d 2}}}{\left(1+{\frac d 4}\right) \left(\left(1+{\frac 4 d}\right) {\overline{\lambda_{j}}}\right) ^{\frac{d}{2}}}
\end{equation}
and
\begin{equation}\label{P10}
{\mathcal N}(z) \ge j  \left(\frac{z}{\left(1+{\frac 4 d}\right) {\overline{\lambda_{j}}}}\right)^{\frac{d}{2}}.
\end{equation}
Moreover, for all $k \ge j \ge 1$,
\begin{equation}\label{beautiful}
\lambda_{k+1}/{\overline{\lambda_{j}}} \le \left(1 + {\frac 4 d}\right) \left({\frac k j}\right)^{\frac 2 d}.
\end{equation}
\end{corollary}

\begin{proof}
The first statement comes from substituting $z = z_j = \left(1 + {\frac 4 d}\right) {\overline{\lambda_{j}}}$
into Eqs. \eqref{P7} and \eqref{P8} and simplifying.  The next
two statements result from substituting the first statement into \eqref{P5.5} and \eqref{P6}.

The final statement \eqref{beautiful} is automatic if
$\lambda_{k+1} \le  \left(1 + {\frac 4 d}\right)
{\overline{\lambda_{j}}}$. Suppose to the contrary that
$\lambda_{k+1} >  \left(1 + {\frac 4 d}\right)
{\overline{\lambda_{j}}}$. As $z$ increases to $\lambda_{k+1}$,
\eqref{P10} becomes valid and in the limit reads
$$
k \ge j  \left(\frac{\lambda_{k+1}}{\left(1+{\frac 4 d}\right) {\overline{\lambda_{j}}}}\right)^{\frac{d}{2}},
$$
which when solved for
$\lambda_{k+1}/{\overline{\lambda_{j}}}$ yields the claim.
\end{proof}

The case $j=k$ reproduces a straightforward and well-known
simplification of Yang's inequality \cite{Yang}, {\it viz.},
\begin{equation}\label{yang2} \lambda_{k+1} \le
\left(1+{\frac 4 d}\right) \, \overline{\lambda}_{k}.
\end{equation}
The case $j=1$ reduces to the Cheng-Yang bound \eqref{CYbound}.
The other cases are new.

We turn now to the relation of Riesz means having
values of $\sigma$ not necessarily differing by integers as in
Theorem~\ref{BasicIneq}. Because the upper and lower bounds in
\eqref{RUprLower} obey the same power law, all Riesz means for
$\sigma \ge 2$ are comparable in the sense that for any pair
$\sigma_{1,2} \ge 2$, and all $z$ larger than a certain value,
$$
c_{\ell} \le \frac{R_{\sigma_1}(z) z^{-\sigma_1}}{R_{\sigma_2}(z) z^{-\sigma_2}} \le c_u,
$$
for explicit non-zero quantities $c_{\ell,u}$ depending only on
$\sigma_{1,2}$, the dimension $d$, and the dimensionless quantity
$|\Omega| \lambda_1^{\frac{d}{2}}$. Moreover, the subsequent
corollary, together with the Laptev-Weidl version of
Berezin-Li-Yau \cite{LapWeidl2} allows the same claim to be made
for $\sigma_{1,2} \ge 1$, with the upper bound holding even for
$\sigma_{2} \ge 0$.

Other constraints on ratios of Riesz means derive from H\"older's inequality,
by which $R_{\sigma}(z)$ is a log-convex function of $\sigma$.  That is,
for $\sigma_0 \le \sigma_1 \le \sigma_2$, if $t$ is chosen so that $\sigma_1 = t \sigma_0 + (1-t) \sigma_2$,
H\"older's inequality applied to
the definition \eqref{Riesz} directly produces:
\begin{equation}\label{logcon}
R_{\sigma_1}(z) \le R_{\sigma_0}(z)^t  R_{\sigma_2}(z)^{1-t}.
\end{equation}
Choosing
$\sigma_0=0, \sigma_1 = \sigma-1$, $\sigma_2 = \sigma$, and $t = \frac 1 \sigma$, we get
$$
R_{\sigma-1}(z) \le \left(R_{\sigma}(z)\right)^{\frac{\sigma-1}{\sigma}} \left(R_0(z)\right)^{\frac 1 \sigma},
$$
or, equivalently,
\begin{equation}\label{Hoelder}
{\mathcal N}(z) \ge \frac{\left(R_{\sigma-1}(z)\right)^{\sigma}}{\left(R_{\sigma}(z)\right)^{\sigma-1}}.
\end{equation}
Along with
 \eqref{Difference2}, we obtain, for
all $z$ and $\sigma \ge 2$,
\begin{equation}
{\mathcal N}(z) \ge \left(\frac{d+ 2 \sigma}{2 \sigma }\right)^{\sigma} z^{-\sigma} R_{\sigma}(z).
\end{equation}
From this stage lower bounds on ${\mathcal N}$, and consequently
upper bounds on $\lambda_{k+1}$, can be derived for all $\sigma
\ge 2$, by optimizing coefficients as before to eliminate
$R_{\sigma}$. The latter, however, are not found to improve upon
the case $\sigma=2$, which corresponds to \eqref{P6} and
\eqref{P10}.

\medskip

\section{Weyl-type estimates of means of eigenvalues} \label{SharpSection}

The Legendre transform is an effective tool for
converting bounds on $R_\sigma(z)$ into bounds on the spectrum,
as has been realized previously, e.g., in
\cite{LapWeidl1}.  We use it in this section to obtain
some improvements on parts of the preceding section, with a
focus on the averages $ \overline{\lambda_k}$.

Recall that if $f(z)$ is a convex function on $\RR^+$ that is superlinear in $z$ as
$z \to +\infty$, its {\it Legendre transform}
\begin{equation}\label{LegDef}
{\mathcal L}\left[f\right](w) := \sup_z{\left\{w z - f(z)\right\}}
\end{equation}
is likewise a superlinear convex function.  Moreover, for each
$w$, the supremum in this formula is attained at some finite value
of $z$.  (A concise treatment of the Legendre transform may be
found for example in chapter 3 of \cite{LCE}.  The Legendre
transform on $\RR^+$ can be understood in terms of the more widely
treated Legendre transform on $\RR$ by restricting to even
functions. ) We also note that $f(z) \ge g(z)$
for all $z \implies {\mathcal
L}\left[f\right](w) \le {\mathcal L}\left[g\right](w)$ for all $w$, and proceed to
calculate the Legendre transform of \eqref{P9}.

First we make a remark about the restriction on the range of
values of $z$ for which Ineq. \eqref{P9} is valid, {\it viz.}, $z
\ge z_j := \left(1 + {\frac 4 d}\right) {\overline{\lambda_{j}}}$.
We can imagine redefining the function on the right for $z \le
z_j$ in some unimportant way, though preserving convexity, so that
it provides a lower bound to the left side for all $z \ge 0$.
While this extended definition contains no information, it ensures
that both sides of  Ineq. \eqref{P5.5} are now defined on a
standard interval guaranteeing that the Legendre transform has the
properties cited above.  Since the maximizing value of $z$ in the
definition \eqref{LegDef} is a nondecreasing function of $w$, it
follows that for $w$ sufficiently large, the maximizing $z$
exceeds $z_j$.  This permits us simply to transform both sides of
Ineq. \eqref{P5.5} to obtain a dual inequality that is valid when
$w \ge w_j$, for some $w_j$, the value of which we shall estimate
{\it post facto}.

The Legendre transforms of the
two sides of Ineq. \eqref{P9} are straightforward:  The maximizing value of $z$ in
the Legendre transform of $R_1$ is attained at one of the critical values,
which, because $R_1$ is piecewise linear, means that
$z_{cr} = \lambda_{k+1}$, and it is easy to check that $k = \left[w\right]$.
Substitution reveals
$$
{\mathcal L}\left[R_1\right](w) = \left(w -  \left[w\right]\right) \lambda_{\left[w\right]+1}
+ \left[w\right] \overline{\lambda_{\left[w\right]}}.
$$

Together with the standard Legendre transform of  the right side  of \eqref{P9}, we obtain
\begin{equation}\label{LegP9a}
\left(w -  \left[w\right]\right) \lambda_{\left[w\right]+1}
+ \left[w\right] \overline{\lambda_{\left[w\right]}}
\le
\frac{2}{j^{\frac 2 d}} \left(\frac{1+\frac d 4}{1+\frac d 2}\right)^{1+\frac 2 d} \overline{\lambda_j} \ w^{1+\frac 2 d}
\end{equation}
for certain values of $w$ and $j$, the determination of which we
now consider. In {Corollary~\ref{lambdabarbds}} it is supposed
that $z \ge \left(1 + {\frac 4 d}\right)
{\overline{\lambda_{j}}}$, and we recall that that there is a
monotonic relationship between $w$ and the maximizing $z_*$ in the
calculation of the Legendre transform of the right side of
\eqref{P9}. (By the maximizing $z_*$ for a Legendre transform we
mean the value such that ${\mathcal L}[f](w) = w z_* - f(z_*)$.)
By an elementary calculation,
\begin{equation}\label{z*}
w = j \left(\frac{1+ \frac d 2}{1+ \frac d 4}\right)\left(\frac{z_*}{\left(1+ \frac 4 d\right) {\overline{\lambda_{j}}}}\right)^{\frac d 2},
\end{equation}
so it follows that if $w$ is restricted to values $\ge j
\frac{1+\frac d 2}{1+\frac d 4}$, then
Ineq.~\eqref{LegP9a} is valid. Meanwhile, for any
$w$ we can always find an integer $k$ such that on the left side
of  \eqref{LegP9a}, $k-1 \le w < k$. If $k > j \frac{1+\frac d
2}{1+\frac d 4}$ and if we let $w$ approach $k$ from below, we
obtain from \eqref{LegP9a}
$$
\lambda_{k}
+ (k-1) \overline{\lambda_{k-1}}
\le
\frac{2}{j^{\frac 2 d}} \left(\frac{1+\frac d 4}{1+\frac d 2}\right)^{1+\frac 2 d} \overline{\lambda_j}k^{1+\frac 2 d}.
$$
The left side of this equation is the sum of the eigenvalues $\lambda_1$ through $\lambda_{k}$, and
the calculation can thus be encapsulated as

\begin{corollary}\label{Leg1}
For $k \ge j \frac{1+\frac d 2}{1+\frac d 4}$, the means of the eigenvalues of the
Dirichlet Laplacian satisfy a universal Weyl-type bound,
\begin{equation}\label{LegP9}
{\overline{\lambda_{k}}}/{\overline{\lambda_{j}}} \le
2 \left(\frac{1+\frac d 4}{1+\frac d 2}\right)^{1+\frac 2 d}\left( {\frac k j}\right)^{\frac 2 d}.
\end{equation}
\end{corollary}

Note that for small values of $j$ the restriction on $k$ is fulfilled simply for $k > j$.  We are aware of only one previous bound somewhat comparable to {Corollary~\ref{Leg1}}, namely in the case
$k=d+1, j=1$, for which Ashbaugh and Benguria \cite{AB93} have shown that
$$
{\overline{\lambda_{d+1}}}/\lambda_1  \le
\frac{d+5}{d+1}.
$$
This is stronger than \eqref{LegP9} when $k=d+1, j=1$. Hence we can set
$j=d+1$ and multiply the two bounds, to obtain:

\begin{corollary}\label{ABHH}
For $k \ge \, \frac{(d+1) \,
(1+\frac{d}{2})}{1+\frac{d}{4}}$,
\begin{equation}\label{P11}
{\overline{\lambda_{k}}}/\lambda_1 \le
\frac{d+5}{2^{\frac 2 d}}\left(\frac{(d+4)}{(d+1)(d+2)}\right)^{1+ \frac 2 d} k^{\frac 2 d}.
\end{equation}
\end{corollary}

We remark that together with \eqref{beautiful}, setting $j=k$, this implies
\begin{equation}\label{P12}
\lambda_{k+1}/\lambda_1 \le
\frac{(d+4)^{2+ \frac 2 d}(d+5)}{2^{\frac 2 d}d(d+1)^{1+ \frac 2 d} (d+2)^{1+ \frac 2 d}} \ k^{\frac 2 d}.
\end{equation}

Finally, we note that any bound in terms of
$\overline{\lambda_k}$, whether from above or below, implies a
bound, with $\overline{\lambda_k}$ replaced by the square root of
$\overline{\lambda_k^2}$ with an additional factor depending only
on dimension, since
\begin{equation}\label{DiscrimBd}
\overline{\lambda_k}^2 \le \overline{\lambda_k^2} \le {\frac{\left({1+ \frac 2 d}\right)^2}{1+ \frac 4 d}} \overline{\lambda_k}^2.
\end{equation}

The first of these inequalities is Cauchy-Schwarz, and the second
is from \cite{HS}, Proposition 6 (ii) (set the parameter $\sigma$
occurring there to 1). The bounds implied by the considerations
above and \eqref{DiscrimBd} for the root-mean-square of the
eigenvalues are of Weyl-type. Moreover, similar statements can be
made about other eigenvalue moments
\[\mathcal{M}_{\sigma}:=\mathcal{M}_{\sigma}(\lambda_1,\lambda_2,\ldots,
\lambda_k):=\left(\frac{1}{k} \, \sum_{\ell \le k}
\lambda_{\ell}^{\sigma}\right)^{1/\sigma}\] when $0<\sigma\le 2$,
as well as for the geometric and harmonic means of the
eigenvalues.
In analogy to the situation of the Riesz means, this is a consequence
of H\"older's inequality
(see also \cite{HLP}), which implies
\begin{equation}
\mathcal{M}_{\sigma}<\mathcal{M}_{\mu} \quad \text{ for }
\sigma<\mu,\notag \end{equation} and the interpolation formula
\begin{equation}
\mathcal{M}_{\sigma}^{\sigma} <
\left(\mathcal{M}_{\mu}^{\mu}\right)^{\frac{\tau-\sigma}{\tau-\mu}}
\,
\left(\mathcal{M}_{\tau}^{\tau}\right)^{\frac{\sigma-\mu}{\tau-\mu}}
\quad \text{ for } \mu<\sigma<\tau. \notag \end{equation}
For $0<\sigma<1$, we specialize to the case
\[\mathcal{M}_{\sigma}<\mathcal{M}_{1}= \overline{\lambda_k}.\]
For $1<\sigma<2$, Ineq. \eqref{DiscrimBd} can be combined with the
interpolation formula for $\mu=1$, $\tau=2$ (Note that
$\mathcal{M}_{2}=\sqrt{\overline{\lambda_k^2}}$.) The claims for
geometric and harmonic means are then clear by virtue of the
classical inequalities relating these means to the arithmetic
mean.

\section{Comparison of universal bounds}\label{comparison}

In this section we compare our universal eigenvalue bounds with
those that have appeared earlier. We consider bounds for
$\overline{\lambda}_{k}$ as well as for $\lambda_{k+1}$. Setting
$j=1$ in \eqref{LegP9}, for $k\ge 2$, we obtain
\begin{equation}\label{simpleP9}
\frac{\overline{\lambda}_{k}}{\lambda_{1}} \le 2
\left(\frac{1+{\frac d 4}}{1+{\frac d 2}}\right)^{1+{\frac 2 d}}
\, k^{\frac 2 d}.
\end{equation}

We claim that \eqref{simpleP9} is sharper than \eqref{CYbound}. To
see this, set $j=0,1, \ldots k-1$ in \eqref{CYbound} and sum, to
get
\begin{equation}
 \sum_{j=0}^{k-1} \frac{\lambda_{j+1}}{\lambda_1} \le
\left(1+{\frac 4 d}\right)  \,
\sum_{j=0}^{k-1} j^{\frac 2 d} \\
\le \left(1+{\frac 4 d}\right) \, \frac{k^{1+{\frac 2
d}}}{1+{\frac 2 d}}, \notag
\end{equation}
where we have simplified the expression by
regarding the middle sum as a left Riemann sum, following
\cite{Her}. After simplification, the Cheng-Yang bound
\eqref{CYbound} implies
\begin{equation}\label{CYbound-av}
\frac{\overline{\lambda}_{k}}{\lambda_1} \le \frac{d+4}{d+2} \,
k^{\frac 2 d}.
\end{equation}
The coefficients of \eqref{CYbound-av} and \eqref{simpleP9} are
plotted against the dimension in Fig. \ref{fig1}. Clearly,
\[
1< 2 \left(\frac{1+{\frac d 4}}{1+{\frac d 2}}\right)^{1+{\frac 2
d}}< \frac{d+4}{d+2}.\]

\begin{figure}[htb!]
  \begin{center}
     \epsfxsize5truein\epsfbox{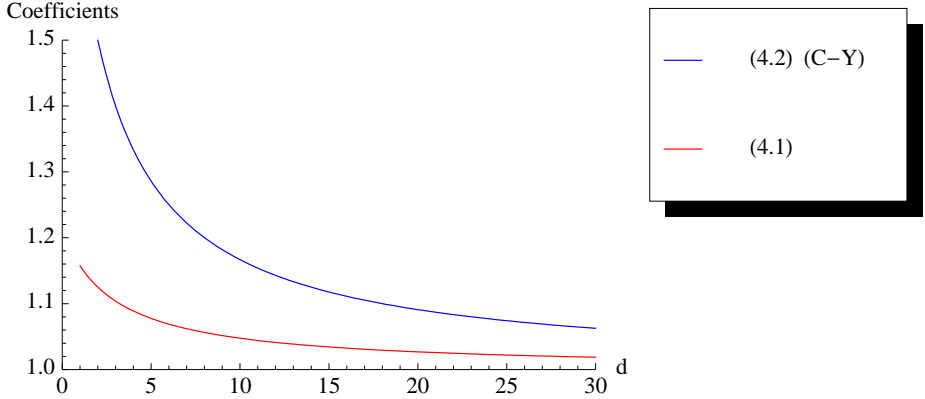}
    \caption{\eqref{simpleP9} vs \eqref{CYbound-av} as a function of $d$ .}
    \label{fig1}
  \end{center}
\end{figure}

The averaged bounds are also plotted in Fig. \ref{fig2} against
$k$ (for $d=4$) along with the bound of Hermi \eqref{her2} and
with \eqref{CYbound-av} and \eqref{simpleP9}, and are seen to
improve the earlier results. The averaged versions of the
Ashbaugh-Benguria bound \eqref{ab94} and that of the Weyl
expression in \eqref{FKbound} are also illustrated along with
\eqref{CYbound-av} and \eqref{simpleP9} in Fig. \ref{fig2}. They
are obtained by dividing the expressions in those equations by
$1+{\frac 2 d}$ (i.e., performing the integration in $k$). This
comparison between the various bounds is displayed numerically for
$\overline{\lambda}_{127}/\lambda_1$ for various dimensions in
Table~\ref{tab1}.

\begin{table}
  \begin{center}
\caption{Bound for $\dfrac{\overline{\lambda}_{127}}{\lambda_1}$
as a function of the dimension $d$.}\label{tab1}
   \begin{tabular}{|c|c|c|c|c|r|}
    \hline \hline
$d$ &   \eqref{simpleP9}  &  \eqref{CYbound-av} & \eqref{her2} & $\overline{AB}$ & $\overline{Weyl}$  \\
    \hline
2 &  142.875   &  190.5   & 163.962 & 339.852 & 43.9204 \\
3 &  27.8886   &  35.3723 & 32.5332 & 89.974  & 8.9804 \\
4 &  12.2686   &  15.0259 & 14.7695 & 40.2459 & 4.0937 \\
5 &  7.48017   &  8.92619 & 9.34082 & 23.3009 & 2.56781 \\
6 &  5.37202   &  6.28316 & 6.95603 & 15.646  & 1.88786 \\
7 &  4.23768   &  4.87795 & 5.67474 & 11.5391 & 1.51906 \\
    \hline
\end{tabular}
\end{center}
\end{table}

\begin{figure}[htb!]
  \begin{center}
    \epsfxsize5truein\epsfbox{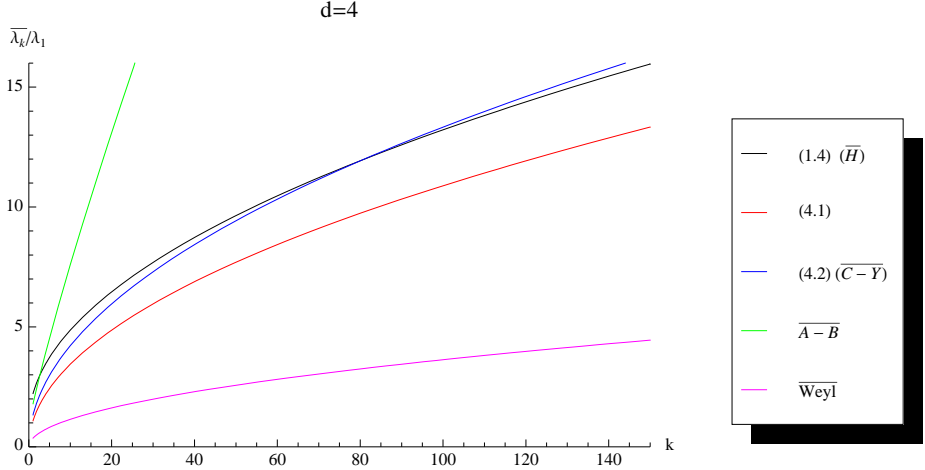}
    \caption{Bound for $\overline{\lambda}_{k}/\lambda_1$ as a function of $k$ ($d=4$).}
  \label{fig2}
  \end{center}
\end{figure}
Notice again the competition between the various bounds with
\eqref{simpleP9} always providing the best constant. One can
also combine \eqref{simpleP9} and \eqref{yang2}. The bounds that
ensue, however, were found not to fare better than \eqref{CYbound}. Note that
we improve by almost 25 percent going from \eqref{CYbound} to
\eqref{simpleP9} when $d=2$; see Table~\ref{table3}.}

\begin{table}
  \begin{center}
\caption{Comparison for coefficients \eqref{simpleP9} and
\eqref{CYbound-av} as a function of the dimension $d$.}
\label{table3}
   \begin{tabular}{|c|c|r|}
    \hline \hline
$d$ & $(4.1)/(\overline{Weyl})$ &  $(4.2)/(\overline{Weyl})$  \\
    \hline
2 &  3.250304  & 4.33739\\
3 &  3.10528  & 3.93884\\
4 &  2.99694  & 3.67049\\
5 &  2.91306  & 3.47619\\
6 &  2.84556  & 3.32818\\
7 &  2.78967  & 3.21116\\
    \hline
\end{tabular}
\end{center}
\end{table}
Similar comparisons can be made for $\overline{\lambda}_{k}$ with
$k\ge d+1$, setting \eqref{P11} beside the bound \eqref{CYbound2}
found in \cite{ChengYang2}.  Again, on average, \eqref{P11} is an
improvement. The comparison is displayed in Table~\ref{table2}
when $k=\Big[\frac{(d+1) \,
(1+\frac{d}{2})}{1+\frac{d}{4}}\Big]+1$. The columns in Table
\ref{table2} represent the ratios between the expressions found in
\eqref{P11} and \eqref{CYbound2} and the averaged version of
\eqref{FKbound}.
\begin{table}
  \begin{center}
\caption{Bound for \eqref{P11} and \eqref{CYbound2} as a function
of the dimension $d$, for $k=\Big[\frac{(d+1) \,
(1+\frac{d}{2})}{1+\frac{d}{4}}\Big ]+1$.}
\label{table2}
   \begin{tabular}{|c|c|r|}
    \hline \hline
$d$ & $(3.5)/(\overline{Weyl})$  &  $(1.9)/(\overline{Weyl})$  \\
    \hline
2 &  2.53014  & 3.08587\\
3 &  2.46466  & 2.92435\\
4 &  2.41249  & 2.80499\\
5 &  2.37103  & 2.71385\\
6 &  2.33756  & 2.64210\\
7 &  2.31003  & 2.58414\\
    \hline
\end{tabular}
\end{center}
\end{table}

\begin{figure}[htb!]
  \begin{center}
    \epsfxsize5truein\epsfbox{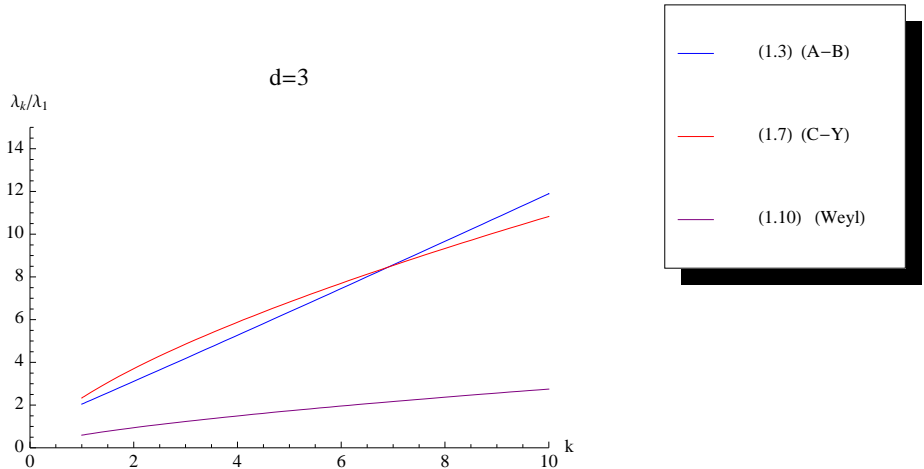}
    \caption{\eqref{ab94} vs \eqref{CYbound} as a function of $k$ ($d=3$).}
    \label{fig3}
  \end{center}
\end{figure}

We conclude by noting, as discussed above, that for low-lying
eigenvalues one cannot expect to improve on the bound \eqref{ab94}
of Ashbaugh-Benguria . This is illustrated in Fig. \ref{fig3}.

\newpage


\end{document}